\documentclass[11pt]{amsart}
\usepackage{geometry}                
\usepackage{graphicx}
\usepackage{amssymb}
\usepackage{epstopdf}
\usepackage{hyperref}
\usepackage{color}
\newcommand{\R}{\mathbb{R}}
\newcommand{\T}{\mathbb{T}}
\newcommand{\C}{\mathbb{C}}
\newcommand{\re}{\operatorname{Re}}

\numberwithin{equation}{section}
\numberwithin{figure}{section}

\theoremstyle{plain} 

\theoremstyle{definition} 

\title{Higher Genus Doubly Periodic Minimal Surfaces}

\author
{Peter Connor}
\address{Peter Connor\\Department of Mathematical Sciences\\Indiana University South Bend\\South
Bend\\IN 46634\\USA\\pconnor@iusb.edu}

\begin{document}
\maketitle

\noindent {\sc Abstract. } {\footnotesize}
We construct Weierstrass data for higher genus embedded doubly periodic minimal surfaces and present numerical evidence that the associated period problem can be solved.  In the orthogonal ends case, there previously was only one known surface for each genus $g$.  We illustrate multiple new examples for each genus $g\geq 3$.  In the parallel ends case, the known examples limit as a foliation of parallel planes with nodes.  We construct a new example for each genus $g\geq 3$ that limit as $g-1$ singly periodic Scherk surfaces glued between two doubly periodic Scherk surfaces and also as a singly periodic surface with four vertical and $2g$ horizontal Scherk ends.


\section{Introduction}
In this paper, we explore the space of embedded doubly periodic minimal surfaces.  Doubly periodic minimal surfaces have top and bottom vertical ends asymptotic to parallel half-planes.  The theory of such surfaces splits into two cases, either the top ends are parallel to the bottom ends or not.  In the non-parallel case, there is classification for genus 0 examples \cite{lm1} and construction methods that produce one surface with orthogonal ends for each genus \cite{ka4}, \cite{ww4}.  In the parallel case, there is classification for up to genus 1 \cite{prt1} and construction methods that prove the existence of many examples of arbitrary genus that limit in foliations by parallel planes with nodes \cite{cw1}.  There are a few other isolated examples with lower genus and parallel ends \cite{rtw1}, \cite{wei2}, \cite{ram1}.

\begin{figure}[t]
	\centerline{ 
		\includegraphics[height=3in]{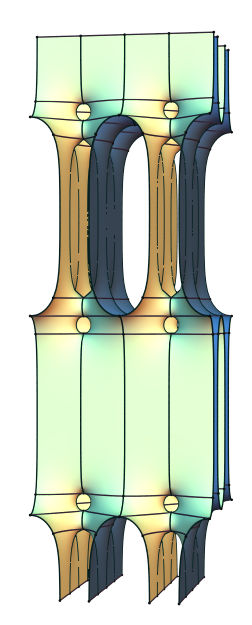}
		\hspace{.25in}
		\includegraphics[height=3in]{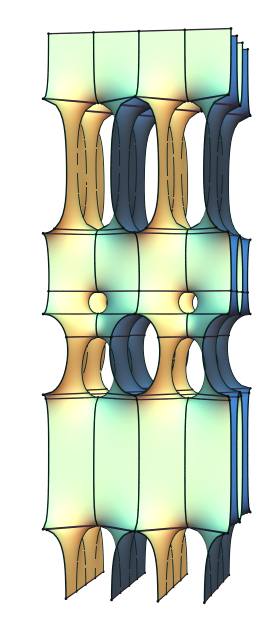}
	}
	\caption{Two new genus 5 examples with orthogonal ends.}
	\label{figure:g5Scherk}
\end{figure}

This paper explores what else is possible beyond the known examples.  Our experiments indicate the existence of new examples with orthogonal ends for genus $g\geq 3$.  See figure \ref{figure:g5Scherk}.  For each genus $g=2n-1$, we conjecture there are $n$ distinct embedded examples.  For each genus $g=2n\geq 4$, we conjecture there are $n-1$ distinct embedded examples.  In the parallel case, we discovered examples of arbitrary genus that don't limit in foliations by parallel planes.

We use the Weierstrass representation for minimal surfaces to develop numerical evidence of the existence of these different types of higher genus embedded doubly periodic minimal surfaces.  The main obstacle to proving the existence of these surfaces is solving the period problem.  We demonstrate how to numerically solve the period problem for each corresponding surface.  The method depends on the surfaces having two orthogonal symmetry planes.  In that case, each surface has a simply connected fundamental domain, from which one can recover the entire surface.  As the period problem involves integrating one-forms over closed curves on a Riemann surface, having a simply connected domain greatly simplifies the calculations. 

The paper is organized as follows.  In section one, we discuss the history of embedded doubly periodic minimal surfaces.

In section two, we outline the Weierstrass representation for doubly periodic minimal surfaces with orthogonal or parallel ends.  Assuming the surfaces have two vertical symmetry planes together with a horizontal symmetry, a template for such surfaces is developed.

In section three, we discuss new higher genus examples with orthogonal ends.  The examples of arbitrary genus constructed by Weber and Wolf in \cite{ww4} have handles between every other pair of ends.  We construct examples in which there are handles between every pair of ends.  A Scherk surface of type $(m,n)$ is a surface with non-parallel Scherk ends that exhibits the alternating pattern of $m$ handles between one pair of ends and $n$ handles between the next pair of ends, and so on.  The Weber-Wolf examples are of type $(0,n)$.  Our experiments indicate the existence of surfaces of type $(1,2n)$ for any $n\geq1$ and of type $(m,n)$ with $m\geq 2$ and $m<n$.  There are, for example, four distinct genus seven surfaces of type $(0,7), (1,6), (2,5)$ and $(3,4)$. As with all doubly periodic minimal surfaces with non-parallel ends, these surfaces lie in a one-parameter family.  The parameter is the angle between the top and bottom ends.  It is possible each genus $g$ example lie in a family that limits to a translation-invariant genus $g$ helicoid, providing evidence that, for each odd genus $g\geq3$ and even genus $g\geq6$, there is more than one translation-invariant helicoid.

In section four, we construct higher genus examples with non-parallel ends with different limit behavior, inspired by the idea of adding multiple handles to a genus one surface with parallel ends.  In particular, these examples add handles to the genus 2 surface from \cite{rtw1}.  The handles all lie along the same vertical symmetry plane, alternating in their direction.  Near one limit, a genus $g$ version looks like $g-1$ singly periodic Scherk surfaces glued in between two doubly periodic Scherk surfaces.  It also limits as a singly periodic surface with four vertical and $2g$ horizontal Scherk ends.  Our experiments suggest that surfaces of this type exist for every genus $g>1$.  We found numerical solutions to the corresponding period problem for genus as high as $g=15$, with no apparent obstructions for higher genus other than dealing with solving a system of large numbers of equations and variables. 

In section five, we discuss possible techniques for proving the existence of the surfaces constructed in sections three and four.

\section{History of embedded doubly periodic minimal surfaces} 
\begin{figure}[h]
	\centerline{ 
		\includegraphics[height=1.4in]{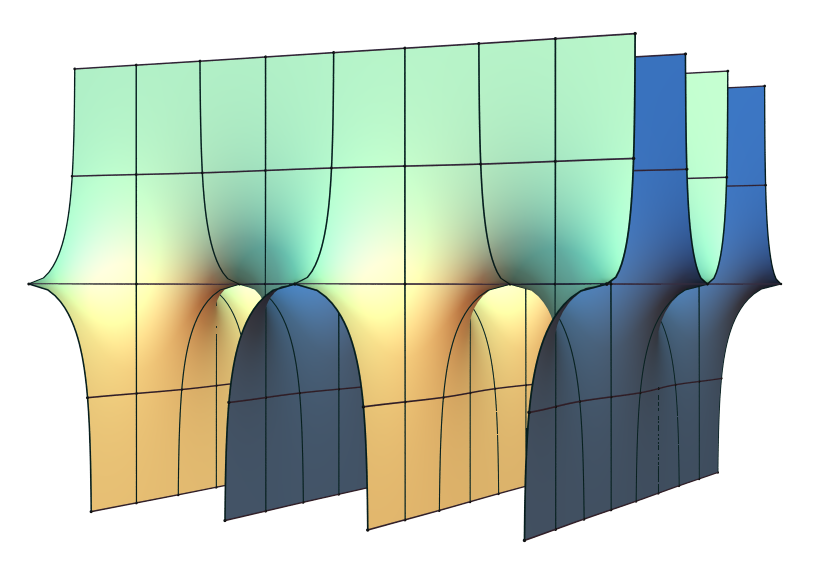}
		\hspace{.5in}
		\includegraphics[height=1.4in]{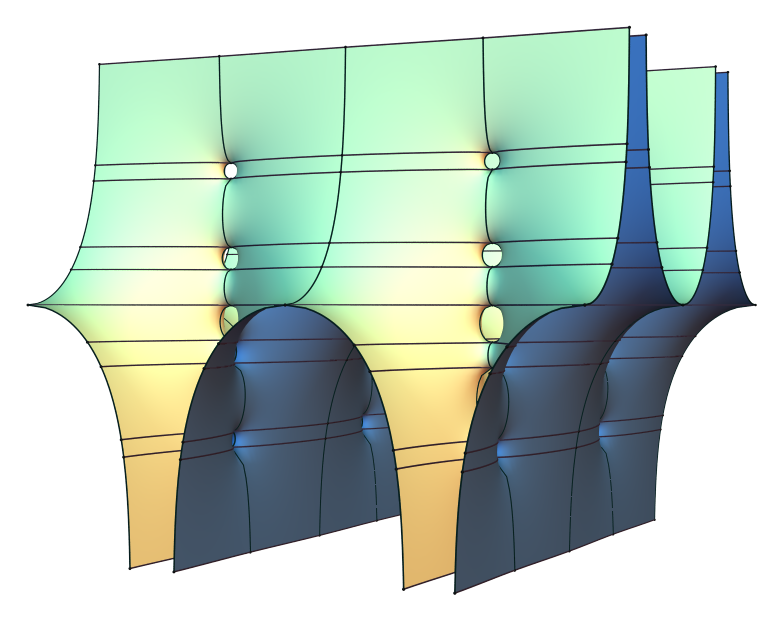}
	}
	\caption{Genus zero and genus four surfaces with orthogonal ends.}
	\label{figure:g0}
\end{figure}

A minimal surface is called {\it doubly periodic} if it is invariant under two linearly independent translations in euclidean space.  After rotating, if necessary, assume the translations are horizontal.  The first example of a doubly periodic minimal surface was discovered in 1835 by Scherk \cite{sche1}, see figure \ref{figure:g0}. 

If $\Lambda$ is the 2-dimensional lattice generated by the two periods of an embedded doubly periodic minimal surface $M$ then its quotient $M/\Lambda$ is an embedded minimal surface in $\R^3/\Lambda$.  Meeks and Rosenberg \cite{mr3} proved that the quotient has a finite number of annular top and bottom ends, referred to as Scherk ends, that are asymptotic to flat annuli.  In $\R^3$ these ends are asymptotic to vertical half-planes.

In order to distinguish between different types of examples, quantities such as genus, number of ends and total curvature for doubly periodic surfaces are measured on the quotient surface.  

Either the top and bottom ends are parallel or non-parallel, with different results and types of examples depending on the case.  As shown by Hauswirth and Traizet \cite{hatr1}, the moduli space of {\it non-degenerate} examples is three-dimensional for parallel ends and one-dimensional for non-parallel ends.  In the parallel case, Meeks and Rosenberg proved that there's an even and equal number of top and bottom ends.

Lazard-Holly and Meeks \cite{lm1} proved that Scherk's doubly periodic surfaces are the only embedded examples of genus 0.  They live in a one-parameter family, with the parameter adjusting the angle $\theta\in(0,\pi/2)$ between the top and bottom ends.  Thus, there are no genus 0 examples with parallel ends.

As the genus increases, the number and complexity of examples increases.  Karcher \cite{ka4} constructed a genus one example with orthogonal ends by adding a handle between every other pair of top (and bottom) ends.  Baginsky and Batista \cite{brb1} and Douglas \cite{dou2} proved that Karcher's example can be deformed by changing the angle $\theta$ between the ends, with $\theta\in(0,\pi/2]$.  As $\theta\rightarrow 0$, the family limits as the translation invariant helicoid with handles.

Weber and Wolf \cite{ww4} generalized Karcher's orthogonal example by constructing one example for each genus $g$ with orthogonal ends and handles between every other pair of ends.  See figure \ref{figure:g0}.  These constitute the only known examples with non-parallel ends, and we demonstrate that there are many more examples by showing how to add handles between each pair of ends.

\begin{figure}[h]
	\centerline{ 
		\includegraphics[height=2in]{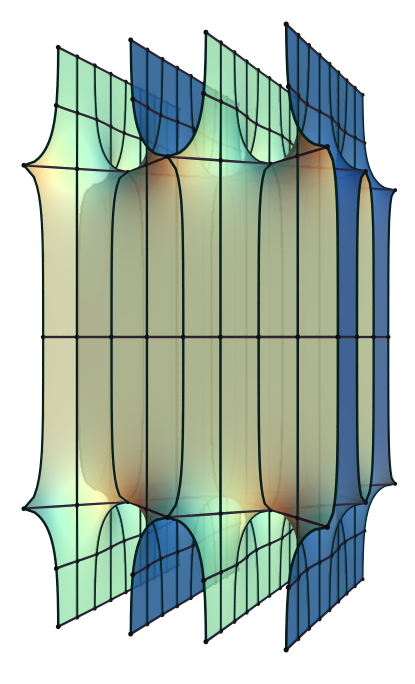}
		\hspace{.5in}
		\includegraphics[height=2in]{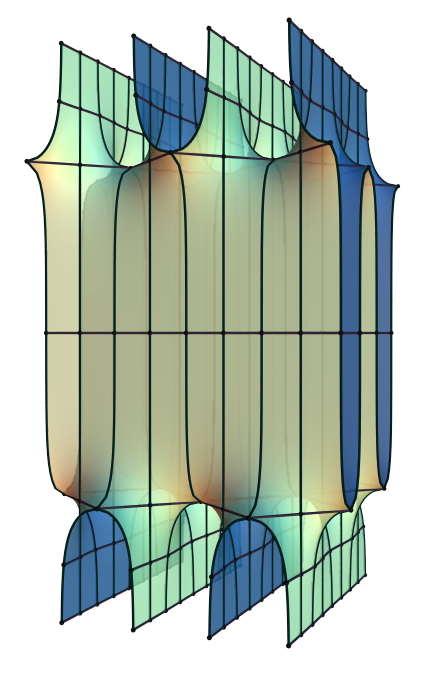}
	}
	\caption{Two examples from the KMR family of surfaces with orthogonal vertical symmetry planes.  The surface on the left also has a horizontal symmetry plane.}
	\label{figure:g1}
\end{figure}

In the parallel case, Karcher \cite{ka4} and Meeks and Rosenberg \cite{mr4} constructed genus one examples with parallel ends.  Perez, Rodriguez, and Traizet \cite{prt1} proved that the three-dimensional moduli space of genus one examples with parallel ends is connected.  These surfaces are referred to as KMR surfaces.  Two KMR surfaces have orthogonal vertical symmetry planes, one of which also has a horizontal symmetry plane.  See figure \ref{figure:g1}. 

\begin{figure}[h]
	\centerline{ 
		\includegraphics[height=2in]{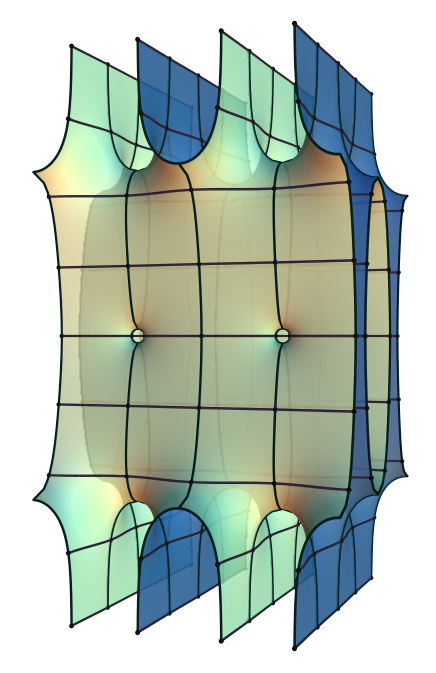}
		\hspace{.5in}
		\includegraphics[height=2in]{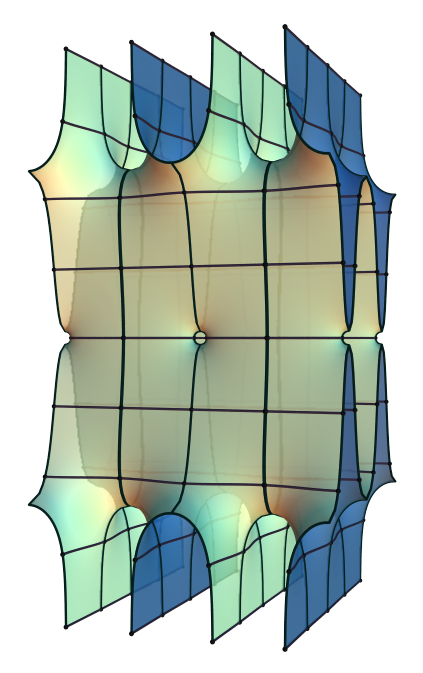}
		\hspace{.5in}
		\includegraphics[height=2in]{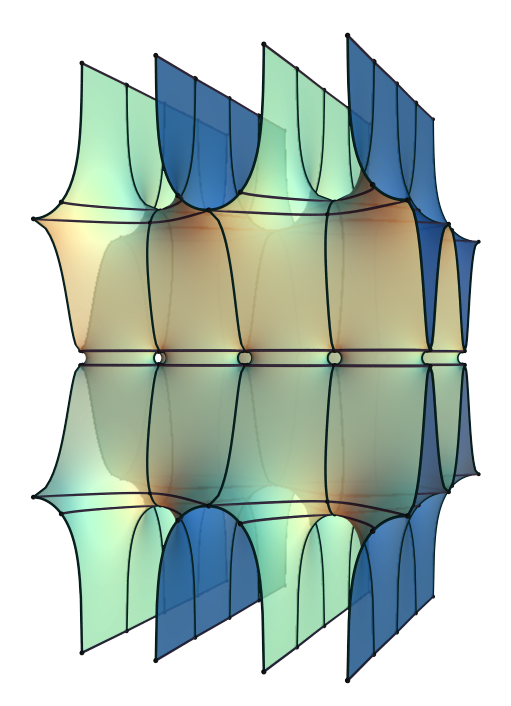}
	}
	\caption{Genus two Wei and RTW surfaces and genus three RTW surface}
	\label{figure:g2}
\end{figure}

Wei \cite{wei2} and Rossman, Thayer, and Wohlgemuth \cite{rtw1} proved the existence of different one-parameter families of genus 2 surfaces  with parallel ends by adding a handle in different ways to a KMR surface.  For convenience, we refer to them as Wei and RTW surfaces.  Additionally, Rossman, Thayer, and Wohlgemuth constructed three different genus three examples, one of which is in figure \ref{figure:g2}.

\begin{figure}[h]
	\centerline{ 
		\includegraphics[height=2in]{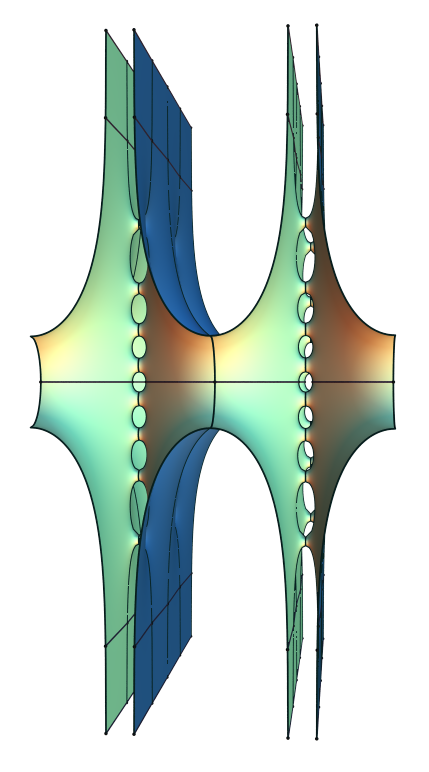}
		\hspace{.5in}
		\includegraphics[height=2in]{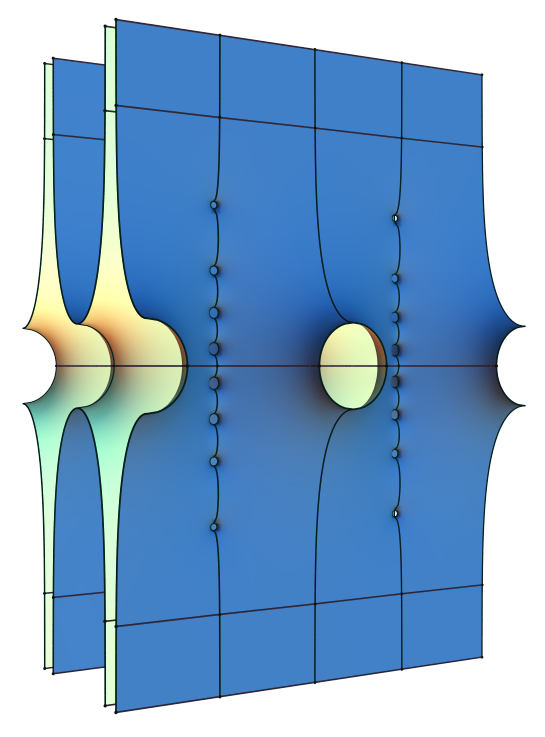}
	}
	\caption{Genus eight surface with parallel ends}
	\label{figure:g8}
\end{figure}

Wei's genus two surface limits as a foliation by vertical parallel planes, and numerical investigations revealed that many higher genus surfaces behaved similarly.  In \cite{cw1}, Connor and Weber proved the existence of new families of embedded, doubly periodic minimal surfaces of arbitrary genus, where each family has a foliation of $\R^3$ by vertical parallel planes as a limit.  In the quotient, these limits can be realized as noded Riemann surfaces, whose components are copies of $\C^*$ with finitely many nodes, and the location of the nodes satisfy a set of balance equations.  Near the limit, the planes move away from each other and are connected by catenoid-shaped necks.  See figure \ref{figure:g8}.

As illustrated in figure \ref{figure:g2}, the Wei and RTW surfaces limit as two doubly periodic Scherk surfaces with a singly periodic Scherk surface glued end-to-end in between.  From another perspective, the genus two Wei and RTW surfaces are constructed by adding a handle to the KMR surface with three perpendicular symmetry planes.  The handles are added along different curves in parallel vertical symmetry planes, with the center of the handles lying in the horizontal symmetry plane.  The genus three RTW surfaces added handles along each of the vertical planar curves.  These constructions led to the idea of constructing higher genus surfaces exhibiting this same behavior by adding more handles along one or both vertical planar curves.

\section{Weierstrass Representation}
Our goal is to construct doubly periodic minimal surfaces with prescribed geometry.  We use the Weierstrass representation for a minimal surfaces.  Let $\tilde{M}$ be a proper minimal surface in $\T\times\R$ with finite genus and four Scherk ends.  By results in \cite{mr3}, $\tilde{M}$ is diffeomorphic to a compact Riemann surface $X$ minus four points, and there exist a meromorphic map $G$ and 1-form $dh$ on $X$ such that the map $f:X\mapsto\T\times\R$ given by
\[
f(z)=\re\int_{z_0}^z\left(\frac{1}{2}\left(\frac{1}{G}-G\right)dh,\frac{i}{2}\left(\frac{1}{G}+G\right)dh,dh\right)
\]
is a minimal immersion whose image is $\tilde{M}$.  The triple $\left(X,G,dh\right)$ is referred to as the Weierstrass representation for $\tilde{M}$.

Conversely, given a Riemann surface $X$, a meromorphic function $G:X\mapsto\overline{\C}$, and a meromorphic one-form $dh$ on $X$, the triple $\left(X,G,dh\right)$ is the Weierstrass representation for a doubly periodic minimal surface with horizontal periods $T_1$ and $T_2$ and vertical Scherk ends if the following conditions hold:
\begin{enumerate}
\item
The zeros of $dh$ are the zeroes and poles of $G$ on $X$ minus the ends, with the same multiplicity.
\item
$dh$ has a poles of order one and $G$ has finite value at the ends.
\item
For each closed curve $\gamma$ on $X$, 
\[
\re\int_{\gamma}\left(\frac{1}{2}\left(\frac{1}{G}-G\right)dh,\frac{i}{2}\left(\frac{1}{G}+G\right)dh,dh\right) =(0,0,0) \;\text{mod}\;\{T_1,T_2\}
\]
This is referred to as the period problem.
\end{enumerate}
We are constructing surfaces with genus and four vertical Scherk ends.  The saddle points of the surface have vertical normal.  Thus, $G$ equals $0$ or $\infty$ at each saddle point.  The top and bottom of each handle are saddle points.  If the genus is $g$ then there will be $2g+2$ saddle points.

Scherk's doubly periodic surface with orthogonal ends can be represented by the Weierstrass data
\[
\begin{split}
X&=\left\{(z,w)\in\overline{\C}^2|w^2=\frac{z-1}{z+1}\right\}\\
G(z)&=w\\
dh&=\frac{1}{z}dz
\end{split}
\]
The saddle points are at $(-1,0)$ and $(1,\infty)$, and the ends are at $(0,\pm i)$ and $(\infty,\pm 1)$.  The bottom ends have normal $(0,\pm 1,0)$, and the top ends have normal $(\pm1,0,0)$.  The Riemann surface $X$ has two automorphisms which correspond to two orthogonal vertical symmetry planes on Scherk's surface.  Thus, a fundamental domain for the surface is given by 
\[
X_1=\left\{(z,w)\in X|z,w>0\right\}
\]

The flat structure of $Gdh$, which is the image of the Schwarz-Christoffel map of $Gdh$ from the domain $X_1$, provides a simple description of the data $(X_1,G,dh)$.  See figure \ref{figure:flatscherk(0,0)}.  The boundaries of the region are the images of the positive and negative real axes.  We can maintain most of the symmetries of Scherk's surface by adding handles along the images of the real axis. On each handle, there are two points $a$ and $b$ with vertical normal vector and $G(a)=0$, $(G(b)=\infty$.  

\begin{figure}[h]
	\centerline{ 
		\includegraphics[height=2in]{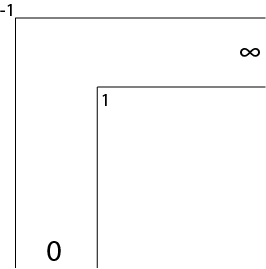}
		\hspace{.25in}
		\includegraphics[height=2in]{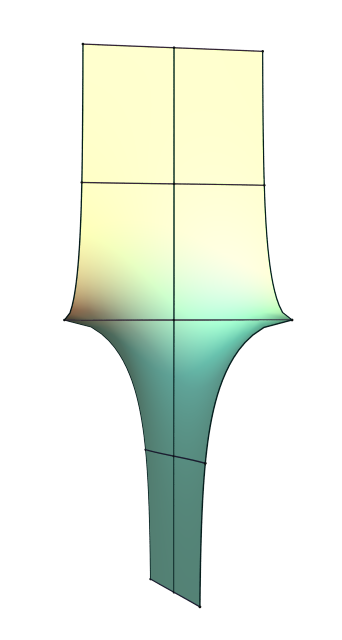}
	}
	\caption{Flat structure for $Gdh$ and image of $X_1$, fundamental domain of Scherk's surface.}
	\label{figure:flatscherk(0,0)}
\end{figure}

A genus $n$ doubly periodic surface can be constructed with Weierstrass data 
\[
\begin{split}
X&=\left\{(z,w)\in\overline{\C}^2|w^2=\prod_{k=1}^{n+1}\frac{z-a_k}{z-b_k}\right\}\\
G(z)&=w\\
dh&=\frac{1}{z}dz
\end{split}
\]
with $a_k,b_k\in\R$ for $k=1,2\ldots, n+1$.  The Riemann surface $X$ has automorphisms $\tau_1(z,w)=(z,-w)$ and $\tau_2(\overline{z},\overline{w})$, and there are the corresponding surface automorphisms
\[
\tau_1^*(\phi_1,\phi_2,\phi_3)=(-\phi_1,\phi_2,\phi_3)
\]
and
\[
\tau_2^*(\phi_1,\phi_2,\phi_3)=(\overline{\phi_1},-\overline{\phi_2},\overline{\phi_3})
\]
of $f(X)$.  Hence, $f(X)$ has two orthogonal vertical symmetry planes.

Note that $G(\infty)=1$, $G(0)=i$ (orthogonal ends) if
\[
b_{n+1}=-a_{n+1}\prod_{k=1}^n\frac{a_k}{b_k}
\]
and $G(0)=1$ (parallel ends) if
\[
b_{n+1}=a_{n+1}\prod_{k=1}^n\frac{a_k}{b_k}
\]

If we add an additional symmetry by requiring that for each $a_k$ there is a $a_j=1/a_k$ or $b_j=1/a_k$ then $X$ has a third automorphism
\[
\tau_3(z,w)=\begin{cases}\left(\frac{1}{z},w\right)\hspace{.1in}\text{if}\hspace{.1in}G(0)=1\\ \left(\frac{1}{z},iw\right)\hspace{.1in}\text{if}\hspace{.1in}G(0)=i\end{cases}
\]
with the corresponding surface automorphism
\[
\tau_3^*(\phi_1,\phi_2,\phi_3)=\begin{cases}-(\phi_1,\phi_2,\phi_3)\hspace{.1in}\text{if}\hspace{.1in}G(0)=1\\ \left(\phi_2,-\phi_1,-\phi_3\right)\hspace{.1in}\text{if}\hspace{.1in}G(0)=i\end{cases}
\]
In both cases, the motivation for working with Weierstrass data with these three symmetries is that the period problem then has far fewer equations.

\subsection{Weber-Wolf examples}
The Weierstrass data for a genus $2n$ Weber-Wolf surface is given by
\[
w^2=\frac{z-1}{z+1}\prod_{k=1}^n\frac{(z-a_{2k-1})(z-a_{2k-1}^{-1})}{(z-a_{2k})(z-a_{2k}^{-1})}
\] 
with $0<a_1<a_2<\cdots<a_{2n}<1$.  See figure \ref{figure:Scherk(0,6)}.

\begin{figure}[t]
	\centerline{ 
		\includegraphics[height=3in]{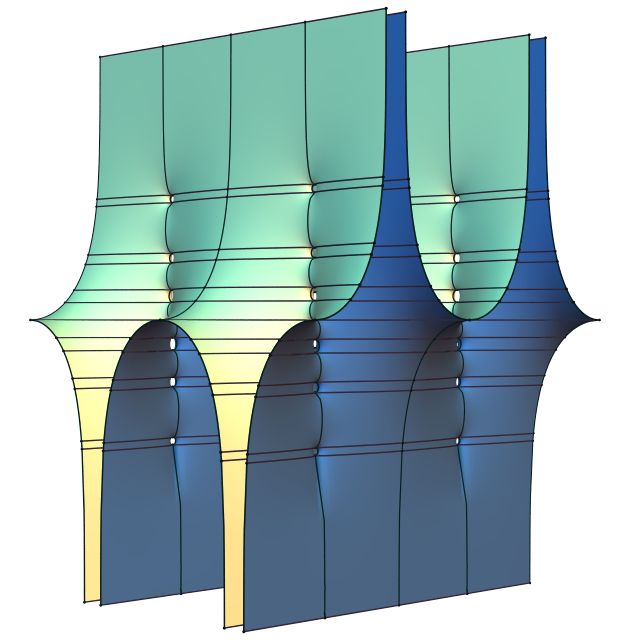}
		\hspace{.25in}
		\includegraphics[height=2.5in]{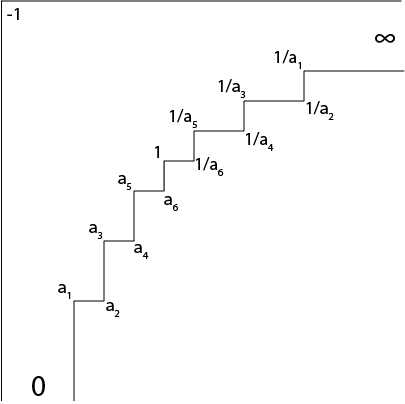}
	}
	\caption{Genus 6 Weber-Wolf surface and the flat structure for $Gdh$.}
	\label{figure:Scherk(0,6)}
\end{figure}

The bottom ends are at $(z,w)=(0,\pm1)$, and the top ends are at $(z,w)=(\infty,\pm i)$.  Note that $G(\infty)=\pm1$ and $G(0)=\pm i$  means the top and bottom ends are orthogonal.  The symmetries of $X$ reduce the period problem to the following equations:
\[
\begin{split}
\int_{a_{2k-1}}^{a_{2k}}\phi_1&=0,\hspace{.1in}k=1,2,\ldots,n\\ 
\int_{a_{2k}}^{a_{2k+1}}\phi_2&=0,\hspace{.1in}k=1,2,\ldots,n-1\\
\int_{a_{2n}}^1\phi_2&=0
\end{split}
\]
The Weierstrass data for a genus $2n+1$ Weber-Wolf surface is given by
\[
w^2=\frac{(z-a_{2n})\left(z-\frac{1}{a_{2n}}\right)}{(z-1)(z+1)}\prod_{k=1}^n\frac{(z-a_{2k-1})(z-a_{2k-1}^{-1})}{(z-a_{2k})(z-a_{2k}^{-1})}
\] 
with $0<a_1<a_2<\cdots<a_{2n+1}<1$.  See figure \ref{figure:Scherk(0,7)}.  The period problem has equations
\[
\begin{split}
\int_{a_{2k-1}}^{a_{2k}}\phi_1&=\int_{a_{2k}}^{a_{2k+1}}\phi_2=0,\hspace{.1in}k=1,2,\ldots,n\\ 
\int_{a_{2n+1}}^1\phi_2&=0
\end{split}
\]
\begin{figure}[h]
	\centerline{ 
		\includegraphics[height=3in]{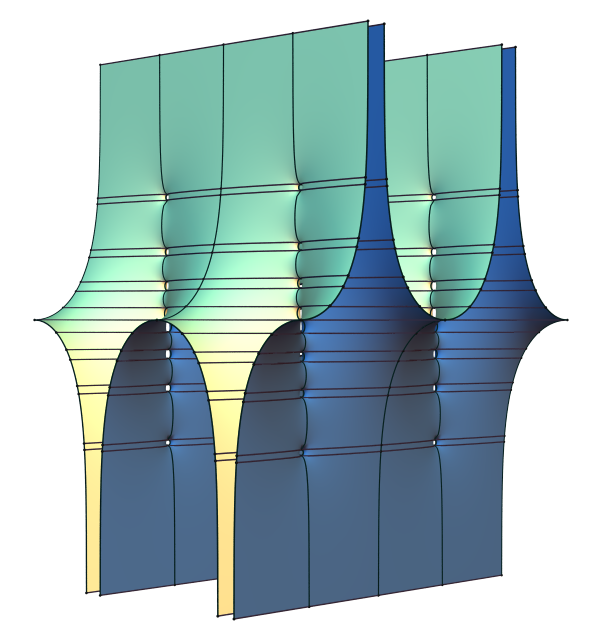}
		\hspace{.25in}
		\includegraphics[height=2.5in]{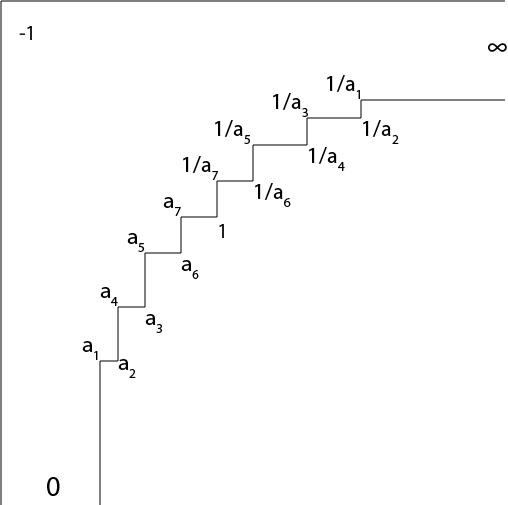}
	}
	\caption{Genus 7 Weber-Wolf surface and the flat structure for $Gdh$.}
	\label{figure:Scherk(0,7)}
\end{figure}

\section{New examples with orthogonal ends}
The Weber-Wolf surfaces add handles between every other pair of ends to Scherk's surface by cutting them out along the image of the positive real axis.  If we also cut them out along the image of the negative real axis then we end up adding handles between every pair of ends to Scherk's surface.  There is numerical evidence that it is possible to do this when the pattern of number of handles between successive pairs of ends is of the form $(1,2n)$.  

\subsection{Surfaces of type $(1,2n)$}
The Weierstrass data of a genus $2n+1$ surface of type $(1,2n)$ is given by
\[
w^2=\frac{(z-1)(z+1)}{(z-b)(z-b^{-1})}\prod_{k=1}^n\frac{(z-a_{2k-1})(z-a_{2k-1}^{-1})}{(z-a_{2k})(z-a_{2k}^{-1})}
\] 
with 
\[
-1<b<0<a_1<a_2<\cdots<a_{2n}<1
\]
The period problem is the same as the Weber-Wolf genus $2n$ surface, with one added period due to the added handle along the negative real axis.
\[
\begin{split}
\int_{a_{2k-1}}^{a_{2k}}\phi_1&=0,\hspace{.1in}k=1,2,\ldots,n\\ 
\int_{a_{2k}}^{a_{2k+1}}\phi_2&=0,\hspace{.1in}k=1,2,\ldots,n-1\\
\int_{a_{2n}}^1\phi_2&=\int_{-1}^b\phi_1=0
\end{split}
\]

\begin{figure}[t]
	\centerline{ 
		\includegraphics[height=3in]{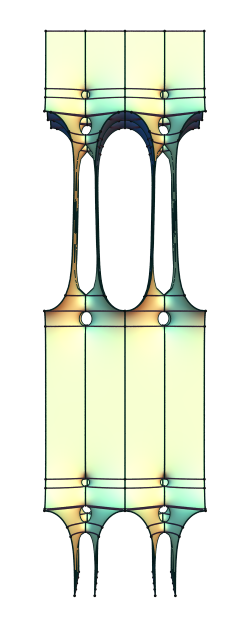}
		\includegraphics[height=3in]{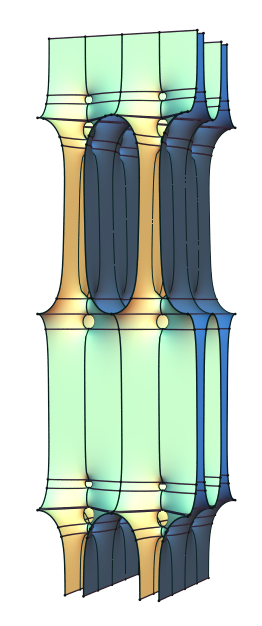}
		\hspace{.25in}
		\includegraphics[height=2.5in]{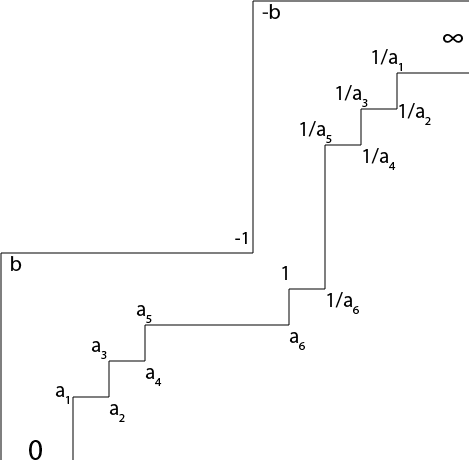}
	}
	\caption{Two views of a Scherk surface of type $(1,6)$ and the flat structure for $Gdh$.}
	\label{figure:Scherk(1,2n+1)}
\end{figure}

This is a system of $2n+1$ equations in $2n+1$ variables which can be solved numerically to find the unique solution.  Note that there is only one solution because these surfaces have non-parallel ends and thus lie in a one-parameter family, with the parameter given by the angle between the top and bottom ends.  We fixed the parameters so the ends are orthogonal.

We used the FindRoot command in Mathematica to find solutions when $n=1,2,3$.  As $n$ increases, it is necessary to increase the working precision.  When $n=3$, using 50 digits of working precision, the period problem is solved when
\[
\begin{split}
a_1&=2.285961038933710244080888845868402093819165097752044509737631058316706\cdot10^{-8}\\
a_2&=4.700161765743336982542669213766119858309585183018657196274538467179463\cdot10^{-8}\\
a_3&=1.8595779281711072657603313530500197096259118694486751650227093422254676\cdot10^{-7}\\
a_4&=8.2857253960766849565137552322547561615808713846473041900491547350816642\cdot10^{-7}\\
a_5&=1.64483174554374477381715594501676269507206942332501684499097006995589289\cdot10^{-6}\\
a_6&=0.29631123476736122914822064611314549156840796261984201940888130050637413082836\\
b&=-1.7178841023135822362953917313799684691657715619679441150357599815555589\cdot10^{-7}\\
\end{split}
\]

\subsection{Surfaces of type $(2,2n)$}
The Weierstrass data of a genus $2n+2$ surface of type $(2,2n)$ is given by
\[
w^2=\frac{(z-1)(z-b_2)\left(z-\frac{1}{b_2}\right)}{(z+1)(z-b_1)\left(z-\frac{1}{b_1}\right)}\prod_{k=1}^n\frac{(z-a_{2k-1})(z-a_{2k-1}^{-1})}{(z-a_{2k})(z-a_{2k}^{-1})}
\] 
with 
\[
-1<b_2<b_1<0<a_1<a_2<\cdots<a_{2n}<1
\]
The period problem is the same as the Weber-Wolf genus $2n$ surface, with two added periods due to the added handles along the negative real axis.
\[
\begin{split}
\int_{a_{2k-1}}^{a_{2k}}\phi_1&=0,\hspace{.1in}k=1,2,\ldots,n\\ 
\int_{a_{2k}}^{a_{2k+1}}\phi_2&=0,\hspace{.1in}k=1,2,\ldots,n-1\\
\int_{a_{2n}}^1\phi_2&=\int_{b_2}^{b_1}\phi_1=\int_{-1}^{b_2}\omega_2=0
\end{split}
\]

\begin{figure}[t]
	\centerline{ 
		\includegraphics[height=3in]{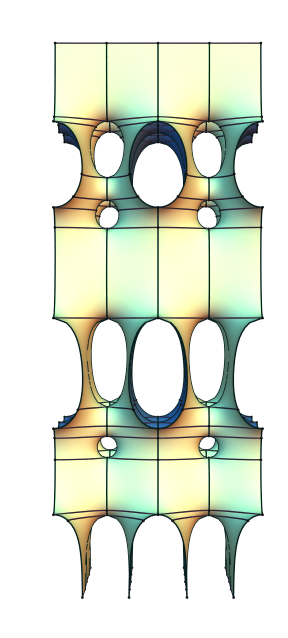}
		\includegraphics[height=3in]{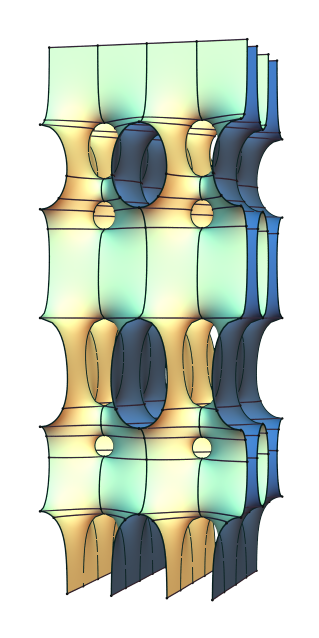}
		\hspace{.25in}
		\includegraphics[height=2.5in]{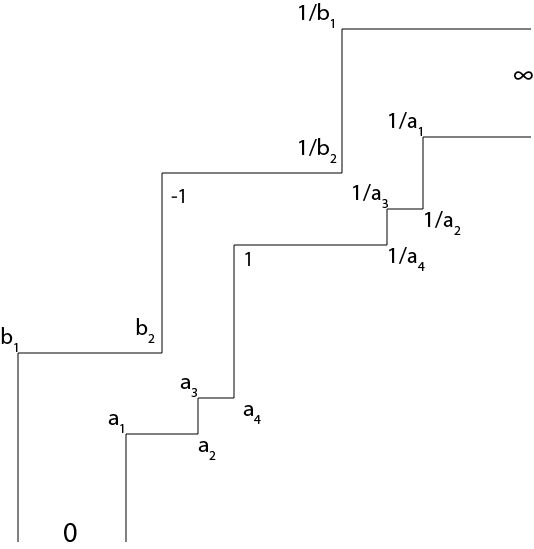}
	}
	\caption{Two views of a Scherk surface of type $(2,4)$ and the flat structure for $Gdh$.}
	\label{figure:Scherk(2,2n)}
\end{figure}

This is a system of $2n+2$ equations in $2n+2$ variables which can be solved numerically to find the unique solution.  We used the FindRoot command with $30$ digits of working precision in Mathematica to find a solution when $n=2$.
\[
\begin{split}
a_1&=5.118542891221798652642273355251580805799073321098417\cdot10^{-6}\\
a_2&=0.0001744687895471917861630587848151868981931348764093613294\\
a_3&=0.0006770540205501970816106707175235603341515662352801700384\\
a_4&=0.004015111358975674034026378269907807153705408517804684585\\
b_1&=-7.363798587282465131067428810155361518052633431720441\cdot10^{-6}\\
b_2&=-0.0012960627225890954277913786548280534111769290970372587867
\end{split}
\]

\subsection{Surfaces of type $(2,2n+1)$}
The Weierstrass data of a genus $2n+3$ surface of type $(2,2n+1)$ is given by
\[
w^2=\frac{(z-a_{2n+1})\left(z-\frac{1}{a_{2n+1}}\right)(z-b_2)\left(z-\frac{1}{b_2}\right)}{(z-1)(z+1)(z-b_1)\left(z-\frac{1}{b_1}\right)}\prod_{k=1}^n\frac{(z-a_{2k-1})(z-a_{2k-1}^{-1})}{(z-a_{2k})(z-a_{2k}^{-1})}
\] 
with 
\[
-1<b_2<b_1<0<a_1<a_2<\cdots<a_{2n+1}<1
\]
The period problem is the same as the Weber-Wolf genus $2n$ surface, with one added period due to the added handle along the negative real axis.
\[
\begin{split}
\int_{a_{2k-1}}^{a_{2k}}\phi_1&=0,\hspace{.1in}k=1,2,\ldots,n\\ 
\int_{a_{2k}}^{a_{2k+1}}\phi_2&=0,\hspace{.1in}k=1,2,\ldots,n\\
\int_{a_{2n}}^1\phi_2&=\int_{-1}^b\phi_1=0
\end{split}
\]

\begin{figure}[t]
	\centerline{ 
		\includegraphics[height=3in]{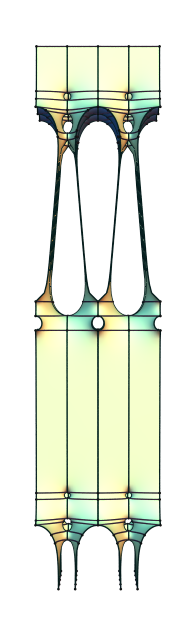}
		\includegraphics[height=3in]{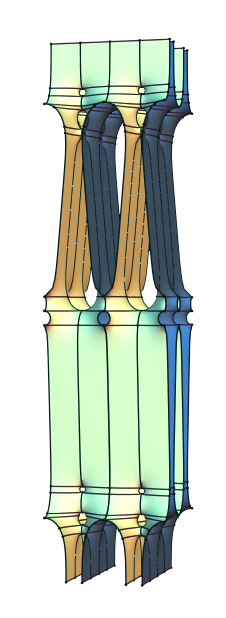}
		\hspace{.25in}
		\includegraphics[height=2.5in]{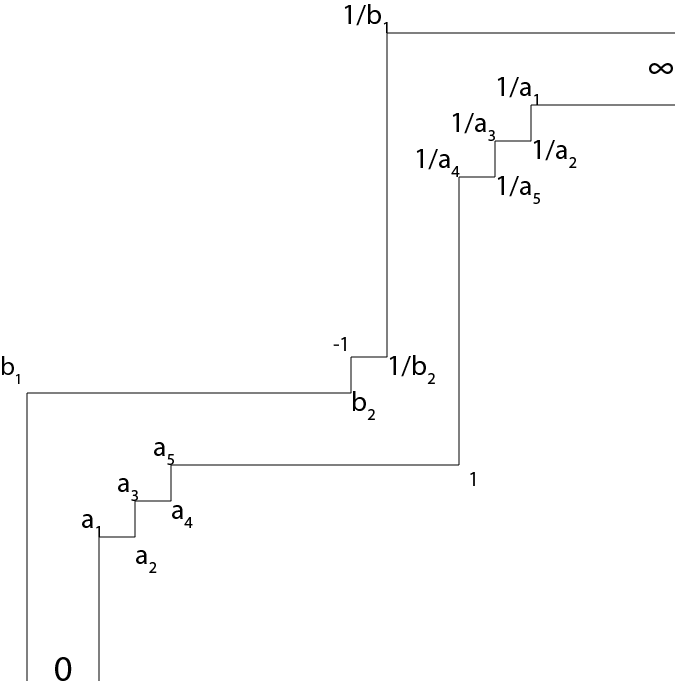}
	}
	\caption{Two views of a Scherk surface of type $(2,5)$ and the flat structure for $Gdh$.}
	\label{figure:Scherk(2,2n+1)}
\end{figure}

We used the FindRoot command in Mathematica to find solutions when $n=1,2$.  When $n=2$, using 40 digits of working precision, the period problem is solved when

\[
\begin{split}
a_1&=1.1282091063550257501145444399176678106500436342667976627322545312066\cdot10^{-10}\\
a_2&=2.5781869298833869544988831073818726093121047370845744173622878810234\cdot10^{-10}\\
a_3&=.900758139358089640789491168379401826694000527832225139272436839872\cdot10^{-10}\\
a_4&=7.52541415527514364659159051880976488113130329298536909943175937427825\cdot10^{-9}\\
a_5&=1.438112750623356054147728184139656280562639183914823161927807185774055\cdot10^{-8}
\end{split}
\]
\[
\begin{split}
b_1&=-6.9941987982384915970848453777676610216381891133676506744108284117738\cdot10^{-10}\\
b_2&=-0.23273985175611771673967941242330200456348758144076705224973842104312950984721
\end{split}
\]

\subsection{Surfaces of type $(3,4)$}
The most exotic example we found was a Scherk surface of type $(3,4)$.  It's Weierstrass data is given by
\[
w^2=\frac{(z-a_{2n+1})\left(z-\frac{1}{a_{2n+1}}\right)(z-b_2)\left(z-\frac{1}{b_2}\right)}{(z-1)(z+1)(z-b_1)\left(z-\frac{1}{b_1}\right)}\prod_{k=1}^n\frac{(z-a_{2k-1})(z-a_{2k-1}^{-1})}{(z-a_{2k})(z-a_{2k}^{-1})}
\] 
with 
\[
-1<b_2<b_1<0<a_1<a_2<\cdots<a_{2n+1}<1
\]
The period problem is the same as the Weber-Wolf genus $2n$ surface, with one added period due to the added handle along the negative real axis.
\[
\begin{split}
\int_{a_{2k-1}}^{a_{2k}}\phi_1&=0,\hspace{.1in}k=1,2\\ 
\int_{a_{2}}^{a_{3}}\phi_2&=\int_{a_{2n}}^1\phi_2=0\\
\int_{b_2}^{b_1}\phi_2&=\int_{b_3}^{b_2}\phi_2=\int_{-1}^{b_3}\phi_1=0
\end{split}
\]

\begin{figure}[h]
	\centerline{ 
		\includegraphics[height=3in]{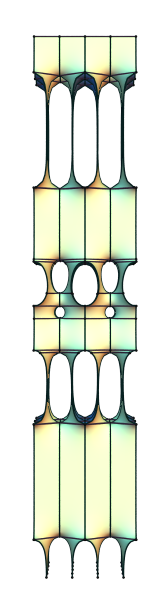}
		\includegraphics[height=3in]{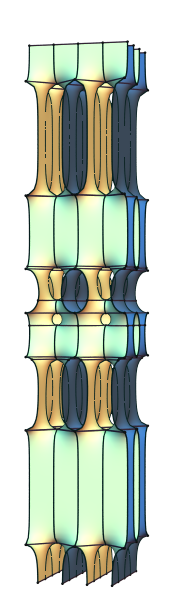}
		\hspace{.25in}
		\includegraphics[height=2.5in]{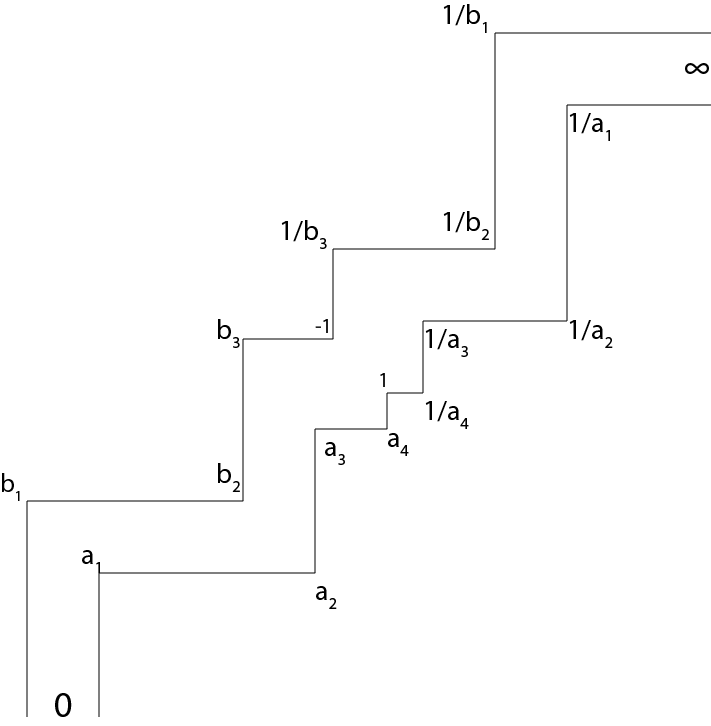}
	}
	\caption{Two views of a Scherk surface of type $(3,4)$ and the flat structure for $Gdh$.}
	\label{figure:Scherk(3,4)}
\end{figure}

This is a system of $7$ equations in $7$ variables which can be solved numerically to find the unique solution.  We used the FindRoot command with 50 digits of working precision in Mathematica to find a solution when 
\[
\begin{split}
a_1&=3.66747800575903677425103971726819961446108179191850037756971821\cdot10^{-13}\\
a_2&=4.89087380310819387522892930330873771443676055616341929003886214\cdot10^{-7}\\
a_3&=0.00318211770512259922630977395737687258793390635072739312898462\\
a_4&=0.21041046622682408387351752969915050930048434243396933164822570\\
b_1&=-3.6763981621447465442422176420782956888909253713440874550960797\cdot10^{-13}\\
b_2&=-5.0663272877810517368694600797246852868580154154927193416664126\cdot10^{-7}\\
b_3&=-0.0043212684921210595889322909764690858329621497375325652336019
\end{split}
\]

\section{New examples with parallel ends}
Near another limit, Wei's surface looks like two doubly periodic Scherk surfaces with a singly periodic Scherk surface glued in between.  Alternatively, one could view it as a KMR surface with a handle.  This led us to consider the possibility that one could glue more than one singly periodic Scherk surface in between two doubly periodic Scherk surfaces.  The building block for these examples are the two KMR surfaces with orthogonal vertical symmetry planes.  Their Weierstrass data is 
\[
w^2=\frac{\left(z+\frac{1}{a}\right)(z-a)}{(z+a)\left(z-\frac{1}{a}\right)}
\]
for the example on the right in figure \ref{figure:g1} or
\[
w^2=\frac{(z-a)\left(z-\frac{1}{a}\right)}{(z+a)\left(z+\frac{1}{a}\right)}
\]
for the example on the left in figure \ref{figure:g1}, with $0<a<1$.

\begin{figure}[t]
	\centerline{ 
		\includegraphics[height=2in]{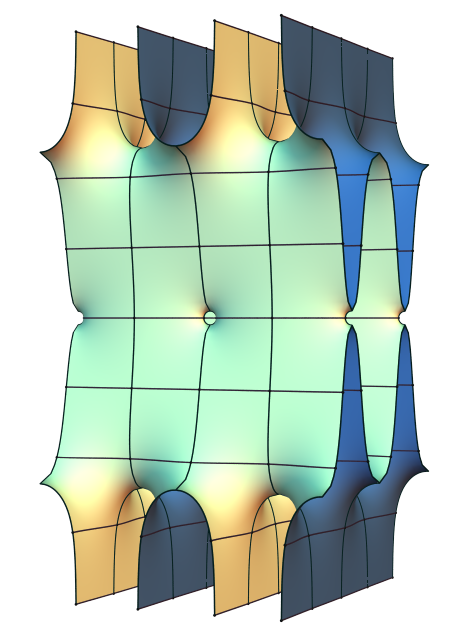}
		\hspace{.25in}
		\includegraphics[height=2in]{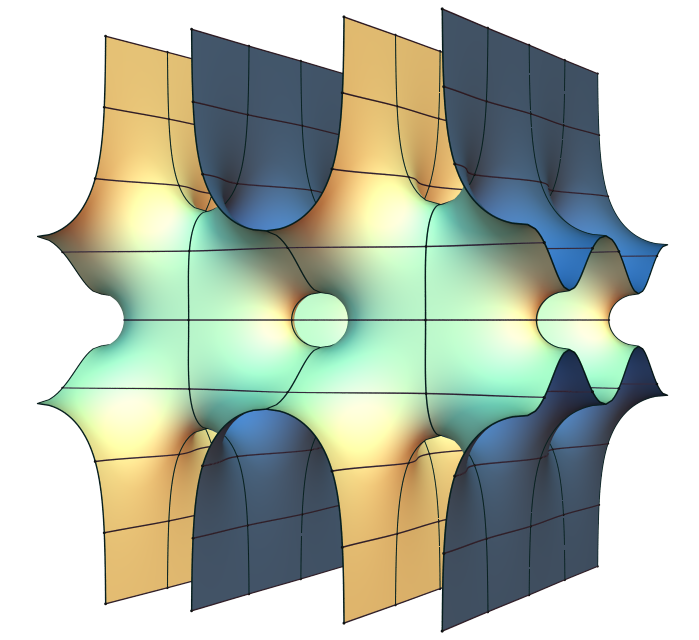}
		\hspace{.25in}
		\includegraphics[height=2in]{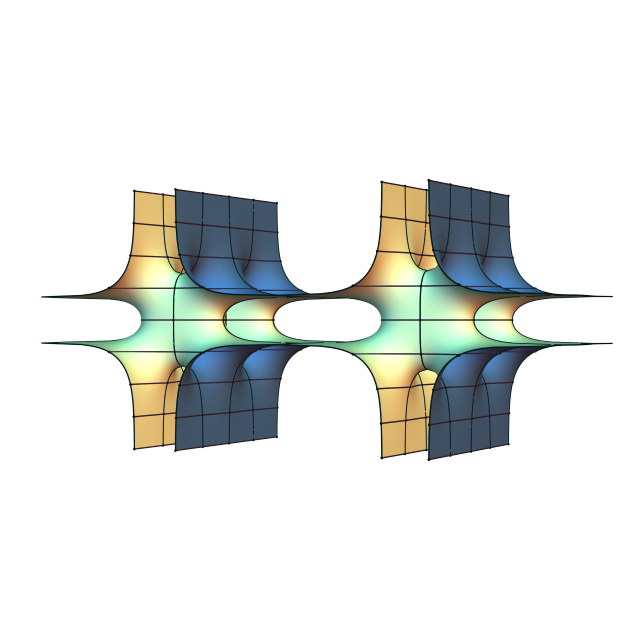}
	}
	\caption{Limits of the RTW family}
	\label{figure:RTWlimits}
\end{figure}

Gluing in extra singly periodic Scherk pieces is realized by simple modifications of $w^2$.  For example, the genus $2$ RTW surface has
\[
w^2=\frac{\left(z-\frac{1}{b}\right)(z-a_1)(z-a_2)}{(z-b)\left(z-\frac{1}{a_1}\right)\left(z-\frac{1}{a_2}\right)}
\]
with $0<b<a_1<a_2<1$.  If $b=-a_1a_2$ then $G(0)=G(\infty)=1$ and the ends are parallel.  Due to the symmetries of the underlying Riemann surface, the period problem for the RTW surface reduces to the equation
\[
\int_{a_1}^{a_2}\phi_2-\pi=0
\]
which has a one-parameter family of solutions, for example
\[
\begin{split}
(a_1,a_2)&=(.0001, 0.7898850561221615)\\
(a_1,a_2)&=(.1, 0.4677900971198217)\\
(a_1,a_2)&=(.265, 0.270905876788826)\\
\end{split}
\]
The RTW surfaces have two discernable limits.  As $(a_1,a_2)\rightarrow(0,1)$, they limit as a singly periodic Scherk surface glued between two doubly periodic Scherk surfaces.  As $(a_1,a_2)\rightarrow\approx(2.6795,2.6795)$, they limit as a singly periodic surface with four vertical and four horizontal Scherk ends.  See figure \ref{figure:RTWlimits}.

\begin{figure}[h]
	\centerline{ 
		\includegraphics[height=3.5in]{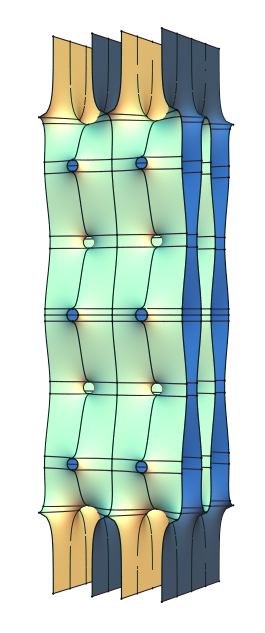}
		\hspace{.05in}
		\includegraphics[height=3.25in]{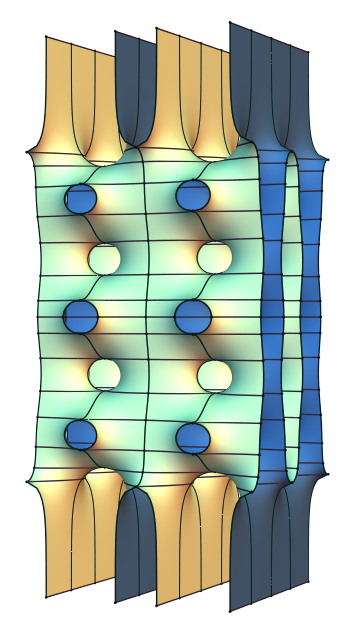}
		\hspace{.05in}
		\includegraphics[height=3in]{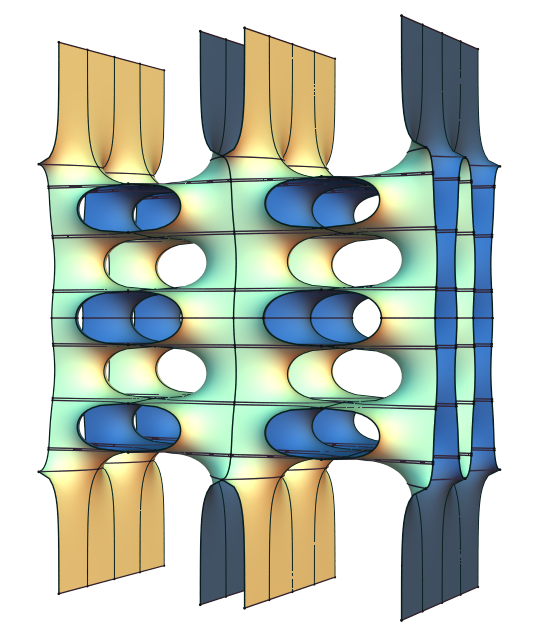}
	}
	\caption{Genus 6 family}
	\label{figure:g6}
\end{figure}

We will add more handles to the RTW surface along the vertical symmetry plane containing its handle, with the handles alternating coming out of and into the surface.  The genus six version, for example, adds four handles above the handle on the RTW surface.  In this case, 
\[
w^2=\frac{\left(z-\frac{1}{b}\right)(z-a_1)(z-a_2)\left(z-\frac{1}{a_3}\right)\left(z-\frac{1}{a_4}\right)(z-a_5)(z-a_6)}{(z-b)\left(z-\frac{1}{a_1}\right)\left(z-\frac{1}{a_2}\right)(z-a_3)(z-a_4)\left(z-\frac{1}{a_5}\right)\left(z-\frac{1}{a_6}\right)}
\]
with $-1<b<0<a_1<a_2<a_3<a_4<a_5<a_6<1$.  If 
\[
b=-\frac{a_1a_2a_5a_6}{a_3a_4}
\]
then $G(0)=G(\infty)=1$.  The period problem reduces to the equations
\[
\begin{split}
\int_{a_1}^{a_2}\phi_2-\pi&=\int_{a_3}^{a_4}\phi_2+\pi=\int_{a_5}^{a_6}\phi_2-\pi=0\\
\int_{a_2}^{a_3}\phi_1&=\int_{a_4}^{a_5}\phi_1=0\\
\end{split}
\]
It has $5$ equations and $6$ variables, so there is a one-parameter family of solutions with the assumed symmetries.  It exhibits similar limiting behavior as the genus 2 RTW family.  One limit is five singly periodic Scherk surfaces glued between two doubly periodic Scherk surfaces.  The other limit is a singly periodic surface with four vertical and twelve horizontal Scherk ends.  See figure \ref{figure:g6}.

\begin{figure}[h]
	\centerline{ 
		\includegraphics[height=4in]{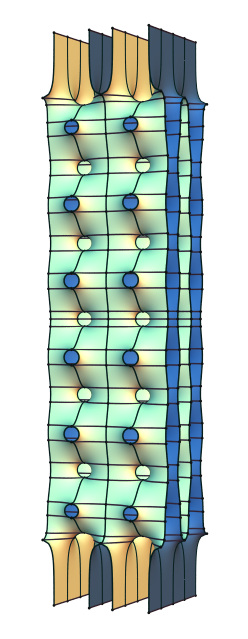}
		\hspace{.05in}
		\includegraphics[height=4in]{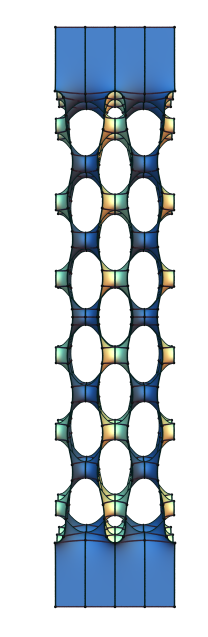}
		\hspace{.25in}
		\includegraphics[height=4in]{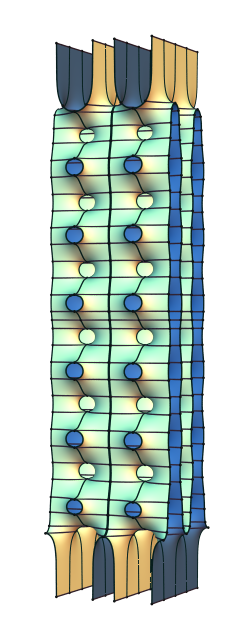}
		\hspace{.05in}
		\includegraphics[height=4in]{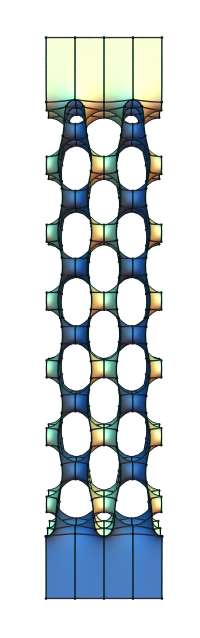}
	}
	\caption{On left: two views of a genus 12 surface.  On right: two views of a genus 13 surface.}
	\label{figure:g12}
\end{figure}

\begin{figure}[h]
	\centerline{ 
		\includegraphics[height=4in]{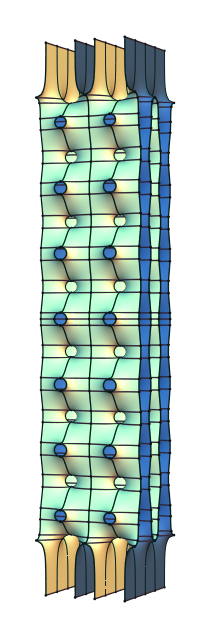}
		\hspace{.05in}
		\includegraphics[height=4in]{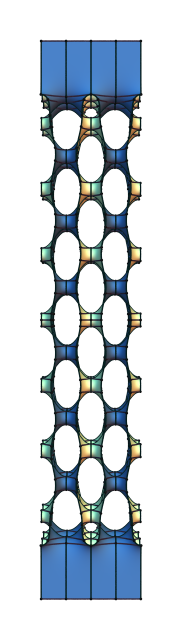}
		\hspace{.25in}
		\includegraphics[height=4in]{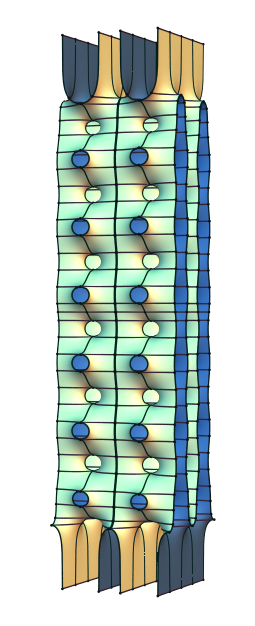}
		\hspace{.05in}
		\includegraphics[height=4in]{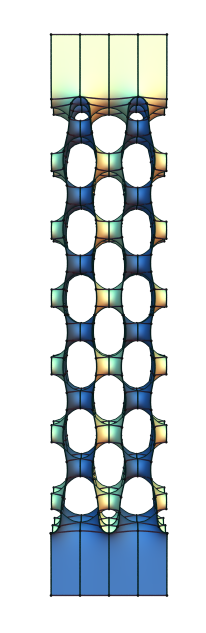}
	}
	\caption{On left: two views of a genus 12 surface.  On right: two views of a genus 13 surface.}
	\label{figure:g14}
\end{figure}

The Weierstrass data for examples with the same characteristics as the above genus six example split into four cases, depending on the genus.  We describe each and demonstrate solutions for surfaces of genus $12, 13, 14$, and $15$.  See figures \ref{figure:g12} and \ref{figure:g14}.

\subsection{Genus=4n}
Let $X$ be the genus $4n$ Riemann surface given by
\[
w^2=\frac{z-\frac{1}{b}}{z-b}\prod_{k=1}^{n}\frac{(z-a_{4k-3})(z-a_{4k-2})\left(z-\frac{1}{a_{4k-1}}\right)\left(z-\frac{1}{a_{4k}}\right)}{\left(z-\frac{1}{a_{4k-3}}\right)\left(z-\frac{1}{a_{4k-2}}\right)(z-a_{4k-1})(z-a_{4k})}
\]
with $-1<b<0<a_1<a_2<\cdots<a_{4n}<1$.  If
\[
b=-\prod_{k=1}^n\frac{a_{4k-3}a_{4k-2}}{a_{4k-1}a_{4k}}
\]
then
$G(\infty)=1$ and $G(0)=1$.

The period problem reduces to the equations
\[
\begin{split}
\int_{a_{2k-1}}^{a_{2k}}\phi_2+(-\pi)^k&=0,\hspace{.1in}k=1,2,\ldots,2n\\ 
\int_{a_{2k}}^{a_{2k+1}}\phi_1&=0,\hspace{.1in}k=1,2,\ldots,2n-1
\end{split}
\]
The period problem has $4n-1$ equations and $4n$ variables.  Thus, we get a one-parameter family of solutions with the assumed symmetries.  This is a subset of the three-dimensional moduli space of these surfaces.  

We used the FindRoot command in Mathematica to find solutions when $n=1,2,3$.  When $n=3$, using 20 digits of working precision, the period problem is solved when
\[
\begin{split}
a_1&=1\cdot10^{-8}\\
a_2&=6.59466639799891609173522939680806465912287592656\cdot10^{-10}\\
a_3&=2.8819141384545408367827643626597393050288222347121\cdot10^{-9}\\
a_4&=3.60735998980815731901024919668500104506308168185473\cdot10^{-8}\\
a_5&=1.878761825242715491218481670214881811480936343264821\cdot10^{-7}\\
a_6&=2.1426402472280116243465380603261391940325586550233778\cdot10^{-6}\\
a_7&=0.0000110243163935562079386000227676184162016976848299812797\\
a_8&=0.0001264421739541800125239658903777500493236340489674605939
\end{split}
\]
\[
\begin{split}
a_9&=0.0006510938026150800869477889312146522995109925417645643327\\
a_{10}&=0.0074648248820113309437057075829391053643379979959846168751\\
a_{11}&=0.0384368739591947293145536401365037438818643908242723861717\\
a_{12}&=0.4406917933970048601093374578213703599711406025352476678762\\
b&=-5.2563\cdot10^{-11}\\
\end{split}
\]

\subsection{Genus=4n+1}
Let $X$ be the genus $4n+1$ Riemann surface given by
\[
w^2=\frac{(z-a_{4n+1})\left(z-\frac{1}{a_{4n+1}}\right)}{(z-b)\left(z-\frac{1}{b}\right)}\prod_{k=1}^{n}\frac{(z-a_{4k-3})(z-a_{4k-2})\left(z-\frac{1}{a_{4k-3}}\right)\left(z-\frac{1}{a_{4k-2}}\right)}{(z-a_{4k-1})(z-a_{4k})\left(z-\frac{1}{a_{4k-1}}\right)\left(z-\frac{1}{a_{4k}}\right)}
\]
with $-1<b<0<a_1<a_2<\cdots<a_{4n+1}<1$, $G(\infty)=1$ and $G(0)=1$.

The period problem reduces to the equations
\[
\begin{split}
\int_{a_{2k-1}}^{a_{2k}}\phi_2+(-\pi)^k&=0,\hspace{.1in}k=1,2,\ldots,2n\\ 
\int_{a_{2k}}^{a_{2k+1}}\phi_1&=0,\hspace{.1in}k=1,2,\ldots,2n-1\\
\int_{1/b_1}^{b_1}\phi_2+\pi&=0
\end{split}
\]
The period problem has $4n$ equations and $4n+1$ variables.  Thus, we get a one-parameter family of solutions with the assumed symmetries.  This is a subset of the three-dimensional moduli space of these surfaces.  

We used the FindRoot command in Mathematica to find solutions when $n=1,2,3$.  As $n$ increases, it is necessary to increase the working precision.  When $n=3$, using 30 digits of working precision, the period problem is solved when
\[
\begin{split}
a_1&=1.2603845439348093259713619930053323710573936489317015427658147213743\cdot10^{-10}\\
a_2&=5.7188763662502975418773254029608668176675750416137360850314932817975\cdot10^{-10}\\
a_3&=3.0065145446521360806465978867857772007627628395052\cdot10^{-9}\\
a_4&=2.09119588124059854086973668043399235755728326827279\cdot10^{-8}\\
a_5&=1.348837690013140059063095058823871025583496552028134\cdot10^{-7}\\
a_6&=8.749514178846944238519225290147618975573247272658584\cdot10^{-7}\\
a_7&=5.5409028958508360611721956233938345329248266678184738\cdot10^{-6}\\
a_8&=0.0000361352083214397670882587515169911541562419402005332588\\
a_9&=0.0002291836892214607884869335205291178203660166494144778651\\
a_{10}&=0.0014939623494188185087754169019070710247635372346728412282\\
a_{11}&=0.0094741073875420897186917809009834149987231949356303425964\\
a_{12}&=0.0617604237670150235055924207337095326652841395655831060664\\
a_{13}&=0.391664016525045117323374978159891090151625490873177707098\\
b&=-6.15\cdot10^{-11}\\
\end{split}
\]

\subsection{Genus=4n+2}
Let $X$ be the genus $4n+2$ Riemann surface given by
\[
w^2=\frac{\left(z-\frac{1}{b}\right)(z-a_{4n+1})(z-a_{4n+2})}{(z-b)\left(z-\frac{1}{a_{4n+1}}\right)\left(z-\frac{1}{a_{4n+2}}\right)}\prod_{k=1}^{n}\frac{(z-a_{4k-3})(z-a_{4k-2})\left(z-\frac{1}{a_{4k-1}}\right)\left(z-\frac{1}{a_{4k}}\right)}{\left(z-\frac{1}{a_{4k-3}}\right)\left(z-\frac{1}{a_{4k-2}}\right)(z-a_{4k-1})(z-a_{4k})}
\]
with $-1<b<0<a_1<a_2<\cdots<a_{4n+2}<1$.  If
\[
b=-a_{4n+1}a_{4n+2}\prod_{k=1}^n\frac{a_{4k-3}a_{4k-2}}{a_{4k-1}a_{4k}}
\]
then
$G(\infty)=1$ and $G(0)=1$.

The period problem reduces to the equations
\[
\begin{split}
\int_{a_{2k-1}}^{a_{2k}}\phi_2+(-\pi)^k&=0,\hspace{.1in}k=1,2,\ldots,2n+1\\ 
\int_{a_{2k}}^{a_{2k+1}}\phi_1&=0,\hspace{.1in}k=1,2,\ldots,2n
\end{split}
\]
The period problem has $4n+1$ equations and $4n+1$ variables.  Thus, we get a one-parameter family of solutions with the assumed symmetries.  This is a subset of the three-dimensional moduli space of these surfaces.  

We used the FindRoot command in Mathematica to find solutions when $n=1,2,3$.  As $n$ increases, it is necessary to increase the working precision.  When $n=3$, using 30 digits of working precision, the period problem is solved when
\[
\begin{split}
a_1&=1\cdot10^{-12}\\
a_2&=7.07635966315915704873463008934708186711541490233514028221304488204\cdot12^{-12}\\
a_3&=2.997319810868618257304188415086578680945894001619918412861378473377\cdot10^{-11}\\
a_4&=4.2213139787159171905867005999520683792521227137531684223382510698251\cdot10^{-10}\\
a_5&=2.1190429251130421258945370091628683438030718065476\cdot10^{-9}\\
a_6&=2.70746442164935563862848528189779616546391459708886\cdot10^{-8}\\
a_7&=1.344092036272544273074732452702322768518584704777767\cdot10^{-7}\\
a_8&=1.7270154143394200845728224769922392130140249814357836\cdot10^{-6}\\
a_{9}&=8.5795767180575296326300138949781821942622640473392949\cdot10^{-6}\\
a_{10}&=0.0001101991510220068176278157452735165345731696998543841768\\
a_{11}&=0.0005474302296481434262712761793407726101973406407808834894\\
a_{12}&=0.0070315475004170885785156441329838010580023202190107441296\\
a_{13}&=0.0349303266765551193911815925592835156330078114457572291206\\
a_{14}&=0.4486670478785153158695747638211060118775881380465661527055\\
b&=-5.3210542942436699788474722194\cdot10^{-13}\\
\end{split}
\]

\subsection{Genus=4n+3}
Let $X$ be the genus $4n+1$ Riemann surface given by
{\footnotesize
\[
w^2=\frac{(z-a_{4n+1})\left(z-\frac{1}{a_{4n+1}}\right)(z-a_{4n+2})\left(z-\frac{1}{a_{4n+2}}\right)}{(z-b)\left(z-\frac{1}{b}\right)(z-a_{4n+3})\left(z-\frac{1}{a_{4n+3}}\right)}\prod_{k=1}^{n}\frac{(z-a_{4k-3})(z-a_{4k-2})\left(z-\frac{1}{a_{4k-3}}\right)\left(z-\frac{1}{a_{4k-2}}\right)}{(z-a_{4k-1})(z-a_{4k})\left(z-\frac{1}{a_{4k-1}}\right)\left(z-\frac{1}{a_{4k}}\right)}
\]}
with $-1<b<0<a_1<a_2<\cdots<a_{4n+3}<1$, $G(\infty)=1$ and $G(0)=1$.

The period problem reduces to the equations
\[
\begin{split}
\int_{a_{2k-1}}^{a_{2k}}\phi_2+(-\pi)^k&=0,\hspace{.1in}k=1,2,\ldots,2n\\ 
\int_{a_{2k}}^{a_{2k+1}}\phi_1&=0,\hspace{.1in}k=1,2,\ldots,2n+1\\
\int_{1/b_1}^{b_1}\phi_2+\pi&=0
\end{split}
\]
The period problem has $4n+2$ equations and $4n+3$ variables.  Thus, we get a one-parameter family of solutions with the assumed symmetries.  This is a subset of the three-dimensional moduli space of these surfaces.  

We used the FindRoot command in Mathematica to find solutions when $n=1,2,3$.  When $n=3$, using 30 digits of working precision, the period problem is solved when
\[
\begin{split}
a_1&=5\cdot10^{-13}\\
a_2&=6.436520914369399010313093473542043713123954908753227224982549961737\cdot10^{-11}\\
a_3&=4.1421218392560424305384662879094426374509776161115104210650197583147\cdot10^{-10}\\
a_4&=1.76403923234234554222669653333257398272696179907630217065265282696619\cdot10^{-9}\\
a_5&=1.42857129482822010518047739096784914147305781814065\cdot10^{-8}\\
a_6&=5.7900815150262660734080873667253934134780659091523\cdot10^{-8}\\
a_7&=4.574845505974320952702888162885640701979832299714309\cdot10^{-7}\\
a_8&=1.8624283891741768799964421591262006611337697726417098\cdot10^{-6}\\
a_9&=0.0000147497857443538071787049578827677906605487618580117471\\
a_{10}&=0.0000600208981480012792979575814087503120656942024970255935\\
a_{11}&=0.0004752444477469941540131623960198059290994726098148511331\\
a_{12}&=0.0019340274806643505816264587782651208086316182766871696561\\
a_{13}&=0.0153117850311046495752669687674177658309850492162375720489\\
a_{14}&=0.0622950460650905819931710051229826193747529336405197641369\\
a_{15}&=0.4939538285674416552069098964089949833351677249500713975416\\
b&=-8.96535437397612068858846648391102561934217329894065034481409536834\cdot10^{-12}\\
\end{split}
\]

\section{Possible methods to prove existence}
There are several different known methods of proving the existence of a minimal surface that may help with proving the existence of the surfaces discussed in sections four and five.  In the orthogonal case, the examples aren't close to a limit of their respective one-parameter families.  Thus, Traizet's regeneration technique won't work here.  The method employed by Weber and Wolf in \cite{ww4} using the flat structures of $Gdh$ and $1/Gdh$ to solve the period problem could work in this setting.  This method uses induction on genus, and so it depends on their being a related example of each genus.  However, the examples with mixed handles don't always have a related example of each genus, for example, there are examples of type $(1,2n)$ but not of type $(1,2n+1)$.  

In the parallel case, Weber and Wolf's handle addition technique from \cite{ww4} is a strong candidate for proving the existence of the surfaces discussed in section five.  We made two unsuccessful attempts to employ Traizet's regeneration technique used in \cite{cw1} to prove the existence of these surfaces in a neighborhood of the limit in which they look like copies of Scherk's singly periodic surface glued in between two copies of Scherk's doubly periodic surface.  One primary difficulty is that, as the angle between the ends of Scherk's singly periodic surfaces goes to zero, the surface converges to a catenoid.  Thus, you go from having four ends to only two ends, losing key data at the limit.  If one can come up with a reasonable expression for the Weierstrass data for a singly periodic surface with four veritcal and $2g$ horizontal Scherk ends then it may work well to apply Traizet's regeneration technique in a neighborhood of that limit. 
\addcontentsline{toc}{section}{References}
\bibliographystyle{plain}
\bibliography{minlit}

\label{sec:liter}

\end{document}